\documentclass[10pt,a4paper]{amsart} 
\usepackage{geometry}
\usepackage{tabls}
\usepackage{url}
\usepackage[utf8]{inputenc}
\usepackage{amsfonts}
\usepackage{amssymb}
\usepackage{mathrsfs}
\usepackage{hyperref}
\hypersetup{hidelinks,
	colorlinks=true,
	allcolors=black,
	pdfstartview=Fit,
	breaklinks=true}
\usepackage{amsmath}
\usepackage{amscd}
\usepackage{fancyhdr}
\usepackage{multirow}
\usepackage{enumerate}
\usepackage{tikz}
\usepackage{tikz-qtree}
\usepackage{paralist}
\usepackage{xypic}
\usepackage{graphicx}
\usepackage{cite}
\usepackage{lipsum}
\usepackage{footmisc}
\usepackage[all,cmtip]{xy}

\allowdisplaybreaks

\numberwithin{equation}{section}
\theoremstyle{plain}
\newtheorem{theorem}{Theorem}[section]
\newtheorem{lemma}[theorem]{Lemma}
\newtheorem{proposition}[theorem]{Proposition}

\newtheorem{remark}[theorem]{Remark}

\numberwithin{equation}{section}

\def\a{\alpha}
\def\b{\beta}
\def\d{\delta}

\def\g{\gamma}

\def\l{\lambda}

\def\n{\vert}

\def\ni{\noindent}

\def\rank{{\rm rank\/}}

\def\R{\mathbb{R} }
\def\C{\mathbb{C}}
\def\P{\mathcal{P}}

\def\N{\mathbb{N}}
\def\Z{\mathbb{Z}}

\def\bz{\bar{z}}
\def\crm{C_{\mathrm{alg}}(\R_{\Theta}^{m})}
\def\incrm{C^\infty_{\mathrm{alg}}(\R_{\Theta}^{m})}
\def\cren{C_{\mathrm{alg}}(\R_{\Theta}^{2n})}
\def\incren{C^\infty_{\mathrm{alg}}(\R_{\Theta}^{2n})}
\def\cron{C_{\mathrm{alg}}(\R_{\Theta}^{2n+1})}
\def\incron{C^\infty_{\mathrm{alg}}(\R_{\Theta}^{2n+1})}
\begin{document}
	
	\title{$K_0$ groups of Connes'  $\Theta$-Deformed $m$-Planes}

	
	%	\institute{
		%		Ren Guan \at
		%		School of Mathematics and Statistics, 
		%		Jiangsu Normal University, Xuzhou 221100, China\\
		%		\email{guanren@jsnu.edu.cn} 
		%	}
	\author{Ren Guan}
	\address{School of Mathematics and Statistics, 
		Jiangsu Normal University, Xuzhou 221100, China}
	\email{guanren@jsnu.edu.cn}

	\begin{abstract}
		We show that the $K_0$ groups of Connes'  $\Theta$-deformed $m$-planes and their smooth versions are all $\Z$.		
	\end{abstract}
	
	
	
	\maketitle
	
	\tableofcontents
	
	\section{Introduction}
	
	In \cite{CD}, Alain Connes and Michel Dubois-Violette propose a new kind of noncommutative  plane, denoted by $C_{\mathrm{alg}}(\R_{\Theta}^{m})$, where $m\in\N^*$ and $\Theta:=\{\theta_{p,q}\}_{1\leq p,q\leq[m/2]}$ a skew-symmetric matrix. For $m=2n$ an even integer, $C_{\mathrm{alg}}(\R_{\Theta}^{2n})$ is the complex unital associative $*$-algebra generated by $2n$ elements $z_p,\bz_q~(p,q=1,2,\ldots,n)$ with relation
	\begin{equation}
		z_pz_q=\lambda_{p,q}z_qz_p,~\bz_p\bz_q=\lambda_{p,q}\bz_q\bz_p,~\bz_pz_q=\lambda_{q,p}z_q\bz_p,~z_p^*=\bz_p
	\end{equation}
	for $p,q=1,2,\ldots,n~(\lambda_{p,q}=e^{i\theta_{p,q}},~\theta_{p,q}=-\theta_{q,p}\in\R)$. For odd $m=2n+1$, $C_{\mathrm{alg}}(\R_{\Theta}^{2n+1})$ is the unital complex $*$-algebra obtained by adding an hermitian generator $x$ to  $C_{\mathrm{alg}}(\R_{\Theta}^{2n})$ with relation $xz_p=z_px,~(p=1,2,\ldots,n)$. $C_{\mathrm{alg}}(\R_{\Theta}^{m})$ are also called the \emph{algebra of complex polynomials on the noncommutative $m$-plane $\R^m_\Theta$}. 
	
	It's not hard to see that every element $T\in\cren$ can be written as a finite sum of the form
	\begin{equation}
		T:=\sum a_{p_1,\ldots,p_{n},q_1,\ldots,q_{n}}z_1^{p_1}\ldots z_n^{p_n}\bz_{1}^{q_{1}}\ldots\bz_{n}^{q_{n}}
	\end{equation}
	with $(p_1,\ldots,p_{n},q_1,\ldots,q_{n})\in\N^{2n}$ and $a_{p_1,\ldots,p_{n},q_1,\ldots,q_{n}}\in\C$. Like the degree of polynomials, we define 
	\begin{equation}
		\deg(z_1^{p_1}\ldots z_n^{p_n}\bz_{1}^{q_{1}}\ldots\bz_{n}^{q_{n}}):=\sum_{k=1}^n(p_r+q_r)
	\end{equation}
	and $\deg(T)$ the maximum  of all degrees of the monomials of $T$. For $\cron$ we have similar representation and definition of degree. When regarded as a subalgebra of the bounded operators $B(H)$ on a separable Hilbert space $H$, mimick the definition of the smooth noncommutative tori\cite{EH}, the smooth version $\incren$ of $\cren$ can be defined as the $C^*$-algebra of the following formal series
	\begin{equation}
		\sum_{p_1,\ldots,p_{n},q_1,\ldots,q_{n}=0}^\infty a_{p_1,\ldots,p_{n},q_1,\ldots,q_{n}}z_1^{p_1}\ldots z_n^{p_n}\bz_{1}^{q_{1}}\ldots\bz_{n}^{q_{n}}
	\end{equation}
	where the coefficient function $\N^{2n}\owns(p_1,\ldots,p_{n},q_1,\ldots,q_{n})\mapsto a_{p_1,\ldots,p_{n},q_1,\ldots,q_{n}}\in\C$ belongs to
	the Schwartz space $\mathcal{S}(\N^{2n})$, i.e., there is a constant $C_r$ for every $r\geq1$ such that 
	\begin{equation}
		\sup_{(p_1,\ldots,p_{n},q_1,\ldots,q_{n})\in\N^{2n}}\left(1+\sum_{k=1}^{n}\left(p_r^2+q_r^2\right)\right)^r a_{p_1,\ldots,p_{n},q_1,\ldots,q_{n}}<C_r.
	\end{equation}
	Similarly we can define $\incron$. 
	
	For a $*$-algebra $A$, a \emph{projector}(or \emph{projection}) $p$ of $A$ is a matrix with entries in $A$ and satisfies $p^2=p=p^*$. We denote by $P(A)$ the set of projectors of $A$. For any two projectors $p,q\in P(A)$, define
	\begin{equation}\label{kkk}
		p+q:=\begin{pmatrix}
			p&0\\
			0&q
		\end{pmatrix},
	\end{equation}
	and call $p,q$ \emph{equivalent}, $p\sim q$, if there is a unitary $u\in M_n(A)$ for some suitable $n\in\N$ such that 
	\begin{equation}
		\begin{pmatrix}
			p&0\\
			0&0
		\end{pmatrix}=u\begin{pmatrix}
			q&0\\
			0&0
		\end{pmatrix}u^*.
	\end{equation}
	Then $P(A)/\sim$  forms a semigroup under the operation \eqref{kkk}. The $K_0$ \emph{group} $K_0(A)$ of $A$ is defined as the Grothendieck group of $P(A)/\sim$. Computing the $K_0$ group of an algebra is a natural question, in general it's not easy, even for the commutative one. 
	
	$K_0$ groups are important for noncommutative geometry, Exel shows that Morita equivalent algebras have isomorphic $K_0$ groups\cite{EX}, and the famous Baum-Connes conjecture\cite{BC} is related to the $K$-theory groups $K_j(C^*_r(G))$ where $j=0,1$ of the reduced $C^*$-algebra $C^*_r(G)$ of a locally compact, Hausdorff and second countable group $G$. See \cite{HLS,MY} for developments in Baum-Connes conjecture and \cite{G,HM,HML,V} for various  examples of calculating $K_0$ groups.
	
	In this paper we compute $K_0(\crm)$ and $K_0(\incrm)$. For $\crm$, we show that $P(\crm)=P(\C)$, which implies $K_0(\crm)=K_0(\C)=\Z$. For the smooth case $\incrm$, there is no such result, but we can prove that any $\mathcal{P}\in P(\incrm)$ is unitary equivalent to 
	\begin{equation}
		\begin{pmatrix}
			I_r&0\\
			0&0
		\end{pmatrix}
	\end{equation}
	for some $r\in\N$, which directly proves $K_0(\incrm)=\Z$. In summary, we have the following theorem.
	\begin{theorem}\label{m}
		$K_0(\crm)=K_0(\incrm)=\Z$ for $\forall m\in\N^*$ and $\Theta$.
	\end{theorem}
	\ni\textbf{Acknowledgements.} 
	This research is partially supported by NSFC grants 12201255.
	
	\section{The nonsmooth cases}
	
	Before the calculation of the $K_0$ groups, we prove the following lemma, the "product rule" of $\cren$:
	\begin{lemma}\label{pr}
		Let $(p_1,\ldots,p_{n},q_1,\ldots,q_{n}),(r_1,\ldots,r_{n},s_1,\ldots,s_{n})\in\N^{2n}$, then
		\begin{equation}
			\begin{aligned}
				&\left(z_1^{p_1}\ldots z_n^{p_n}\bz_{1}^{q_{1}}\ldots\bz_{n}^{q_{n}}\right)\left(z_1^{r_1}\ldots z_n^{r_n}\bz_{1}^{s_{1}}\ldots\bz_{n}^{s_{n}}\right)\\
				=&\left(\prod_{l=1}^n\prod_{k=l+1}^n\l_{k,l}^{p_kr_l+q_ks_l+r_kq_l-q_kr_l}\right)z_1^{p_1+r_1}\ldots z_n^{p_n+r_n}\bz_{1}^{q_{1}+s_1}\ldots\bz_{n}^{q_{n}+s_n}.
			\end{aligned}
		\end{equation}
	\end{lemma}
	\begin{proof}
		For $r,s\in\N$, we have
		\begin{equation}
			z_p^rz_q^s=z_q\lambda_{p,q}^rz_p^rz_q^{s-1}=z_q^2\lambda_{p,q}^{2r}z_p^rz_q^{s-2}=\ldots=\lambda_{p,q}^{rs}z_q^sz_p^r.
		\end{equation}
		Similarly, $\bz_p^r\bz_q^s=\lambda_{p,q}^{rs}\bz_q^s\bz_p^r$ and $\bz_p^rz_q^s=\lambda_{q,p}^{rs}z_q^s\bz_p^r=\lambda_{p,q}^{-rs}z_q^s\bz_p^r$. Hence
		\begin{align*}
			&\left(z_1^{p_1}\ldots z_n^{p_n}\bz_{1}^{q_{1}}\ldots\bz_{n}^{q_{n}}\right)\left(z_1^{r_1}\ldots z_n^{r_n}\bz_{1}^{s_{1}}\ldots\bz_{n}^{s_{n}}\right)\\
			=&\l_{1,1}^{p_{1}r_1}\l_{2,1}^{p_{2}r_1}\ldots\l_{n,1}^{p_{n}r_1}\l_{1,1}^{q_{1}r_1}\l_{1,2}^{q_{2}r_1}\ldots\l_{1,n}^{q_{n}r_1} z_1^{p_1+q_1}z_2^{p_2}\ldots z_n^{p_n}\bz_{1}^{q_{1}}\ldots\bz_{n}^{q_{n}}z_2^{r_2}\ldots z_n^{r_n}\bz_{1}^{s_{1}}\ldots\bz_{n}^{s_{n}}\\
			=&\l_{1,1}^{p_{1}r_1}\l_{2,1}^{p_{2}r_1}\ldots\l_{n,1}^{p_{n}r_1}\l_{1,1}^{q_{1}r_1}\l_{1,2}^{q_{2}r_1}\ldots\l_{1,n}^{q_{n}r_1}\l_{2,2}^{p_{2}r_2}\l_{3,2}^{p_{3}r_2}\ldots\l_{n,2}^{p_{n}r_2}\l_{2,1}^{q_{1}r_2}\l_{2,2}^{q_{2}r_2}\ldots\l_{2,n}^{q_{n}r_2} z_1^{p_1+r_1}z_2^{p_2+r_2}z_3^{p_3}\ldots z_n^{p_n}\\
			&\times\bz_{1}^{q_{1}}\ldots\bz_{n}^{q_{n}}z_3^{r_3}\ldots z_n^{r_n}\bz_{1}^{s_{1}}\ldots\bz_{n}^{s_{n}}\\
			=&\ldots\\
			=&\prod_{l=1}^n\left(\prod_{k=l}^n\l_{k,l}^{p_k}\prod_{k=1}^n\l_{l,k}^{q_k}\right)^{r_l}z_1^{p_1+r_1}\ldots z_n^{p_n+r_n}\bz_{1}^{q_{1}}\ldots\bz_{n}^{q_{n}}\bz_{1}^{s_{1}}\ldots\bz_{n}^{s_{n}}\\
			=&\prod_{l=1}^n\left(\prod_{k=l}^n\l_{k,l}^{p_k}\prod_{k=1}^n\l_{l,k}^{q_k}\right)^{r_l}\left(\prod_{k=1}^n\l_{k,1}^{q_k}\right)^{s_1}z_1^{p_1+r_1}\ldots z_n^{p_n+r_n}\bz_{1}^{q_{1}+s_1}\bz_{2}^{q_{2}}\ldots\bz_{n}^{q_{n}}\bz_{2}^{s_{2}}\ldots\bz_{n}^{s_{n}}\\
			=&\ldots\\
			=&\prod_{l=1}^n\left(\prod_{k=l}^n\l_{k,l}^{p_k}\prod_{k=1}^n\l_{l,k}^{q_k}\right)^{r_l}\prod_{l=1}^n\left(\prod_{k=l}^n\l_{r,s}^{q_k}\right)^{s_l}z_1^{p_1+r_1}\ldots z_n^{p_n+r_n}\bz_{1}^{q_{1}+s_1}\ldots\bz_{n}^{q_{n}+s_n}\\
			=&\left(\prod_{l=1}^n\prod_{k=l+1}^n\l_{k,l}^{p_kr_l+q_ks_l+r_kq_l-q_kr_l}\right)z_1^{p_1+r_1}\ldots z_n^{p_n+r_n}\bz_{1}^{q_{1}+s_1}\ldots\bz_{n}^{q_{n}+s_n}.
		\end{align*}
	\end{proof}
	With the help of Lemma \ref{pr}, we have the following crucial proposition.
	\begin{proposition}\label{p2n}
		If $\mathcal{P}\in P(C_{\mathrm{alg}}(\R_{\Theta}^{2n}))$, then $\P\in P(\C)$.
	\end{proposition}
	\begin{proof}
		Let $$\P=\begin{pmatrix}
			p_{1,1}&\ldots&p_{1,N}\\
			\vdots&\ddots&\vdots\\
			p_{N,1}&\ldots&p_{N,N}\\
		\end{pmatrix}\in M_N(C_{\mathrm{alg}}(\R_{\Theta}^{2n}))$$
		be an $N$-dimensional projector of $C_{\mathrm{alg}}(\R_{\Theta}^{2n})$. Then $\P^2=\P=\P^*$, and
		\begin{align*}
			\begin{pmatrix}
				p_{1,1}&\ldots&p_{1,N}\\
				\vdots&\ddots&\vdots\\
				p_{N,1}&\ldots&p_{N,N}\\
			\end{pmatrix}=\P=\P^2&=\P\P^*\\
			&=\begin{pmatrix}
				p_{1,1}&\ldots&p_{1,N}\\
				\vdots&\ddots&\vdots\\
				p_{N,1}&\ldots&p_{N,N}\\
			\end{pmatrix}\cdot\begin{pmatrix}
				p_{1,1}^*&\ldots&p_{N,1}^*\\
				\vdots&\ddots&\vdots\\
				p_{1,N}^*&\ldots&p_{N,N}^*\\
			\end{pmatrix}\\
			&=\begin{pmatrix}
				\sum_{k=1}^Np_{1,k}p_{1,k}^*&\ldots&\sum_{k=1}^Np_{1,k}p_{N,k}^*\\
				\vdots&\ddots&\vdots\\
				\sum_{k=1}^Np_{N,k}p_{1,k}^*&\ldots&\sum_{k=1}^Np_{N,k}p_{N,k}^*\\
			\end{pmatrix},
		\end{align*}
		so for $1\leq k\leq N$,
		\begin{equation}\label{pll}
			p_{k,k}=\sum_{l=1}^Np_{k,l}p_{k,l}^*.
		\end{equation}
		From Lemma \ref{pr}, for  $(p_1,\ldots,p_{n},q_1,\ldots,q_{n})\in\N^{2n}$, we have
		\begin{align*}
			&\left(z_1^{p_1}\ldots z_n^{p_n}\bz_{1}^{q_{1}}\ldots\bz_{n}^{q_{n}}\right)\left(z_1^{p_1}\ldots z_n^{p_n}\bz_{1}^{q_{1}}\ldots\bz_{n}^{q_{n}}\right)^*\\
			=&z_1^{p_1}\ldots z_n^{p_n}\bz_{1}^{q_{1}}\ldots\bz_{n}^{q_{n}}z_n^{q_n}\ldots z_1^{q_1}\bz_n^{p_n}\ldots\bz_1^{p_1}\\
			=&\left(\prod_{s=1}^{n-1}\prod_{r=s+1}^n\l_{r,s}^{p_rp_s+q_rq_s}\right)z_1^{p_1}\ldots z_n^{p_n}\bz_{1}^{q_{1}}\ldots\bz_{n}^{q_{n}}z_1^{q_1}\ldots z_n^{q_n}\bz_{1}^{p_{1}}\ldots\bz_{n}^{p_{n}}\\
			=&\left(\prod_{s=1}^{n-1}\prod_{r=s+1}^n\l_{r,s}^{p_rp_s+q_rq_s}\right)\prod_{s=1}^n\left(\prod_{r=s}^n\l_{r,s}^{p_rq_s+p_sq_r}\prod_{r=1}^n\l_{s,r}^{q_rq_s}\right)z_1^{p_1+q_1}\ldots z_n^{p_n+q_n}\bz_1^{p_1+q_1}\ldots \bz_n^{p_n+q_n}\\
			=&\left(\prod_{s=1}^{n-1}\prod_{r=s+1}^n\l_{r,s}^{(p_r+q_r)(p_s+q_s)}\right)z_1^{p_1+q_1}\ldots z_n^{p_n+q_n}\bz_1^{p_1+q_1}\ldots \bz_n^{p_n+q_n}.
		\end{align*}
		The above formula implies that for fixed $(M_1,\ldots,M_n)\in\N^{n}$ and any $(p_1,\ldots,p_{n},q_1,\ldots,q_{n})\in\N^{2n}$ such that $p_j+q_j=M_j$ for all $1\leq j\leq n$, we have
		\begin{equation}\label{pqm}
			\left(z_1^{p_1}\ldots z_n^{p_n}\bz_{1}^{q_{1}}\ldots\bz_{n}^{q_{n}}\right)\left(z_1^{p_1}\ldots z_n^{p_n}\bz_{1}^{q_{1}}\ldots\bz_{n}^{q_{n}}\right)^*=\left(\prod_{s=1}^{n-1}\prod_{r=s+1}^n\l_{r,s}^{M_rM_s}\right)z_1^{M_1}\ldots z_n^{M_n}\bz_1^{M_1}\ldots \bz_n^{M_n}.
		\end{equation}
		Let $M:=\max\{\deg (p_{k,l})\}_{1\leq k,l\leq N}$. If $M>0$, set
		\begin{equation}
			p_{k,l}=\sum_{\sum_{j=1}^n(p_j+q_j)\leq M}a_{p_1,\ldots,p_{n},q_1,\ldots,q_{n}}^{k,l}z_1^{p_1}\ldots z_n^{p_n}\bz_{1}^{q_{1}}\ldots\bz_{n}^{q_{n}},
		\end{equation}
		then for $1\leq k\leq N$,
		\begin{align*}
			p_{k,k}=&\sum_{l=1}^Np_{k,l}p_{k,l}^*\\
			=&\sum_{l=1}^N\sum_{\sum_{j=1}^n(p_j+q_j)\leq M}\sum_{\sum_{j=1}^n(p'_j+q'_j)\leq M}a_{p_1,\ldots,p_{n},q_1,\ldots,q_{n}}^{k,l}\overline{a_{p'_1,\ldots,p'_{n},q'_1,\ldots,q'_{n}}^{k,l}}z_1^{p_1}\ldots z_n^{p_n}\bz_{1}^{q_{1}}\ldots\bz_{n}^{q_{n}}\\
			&\times z_n^{q'_n}z_{n-1}^{q'_{n-1}}\ldots z_1^{q'_1}\bz_n^{p'_n}\bz_{n-1}^{p'_{n-1}}\ldots\bz_1^{p'_1}.
		\end{align*}
		Consider the coefficient of $z_1^M\bz_1^M$, we have
		\begin{equation}
			\begin{aligned}
				0=&\sum_{l=1}^N\left(a_{M,0,\ldots,0,0,0,\ldots,0}^{k,l}\overline{a_{M,0,\ldots,0,0,0,\ldots,0}^{k,l}}+a_{0,0,\ldots,0,M,0,\ldots,0}^{k,l}\overline{a_{0,0,\ldots,0,M,0,\ldots,0}^{k,l}}\right)\\
				=&\sum_{l=1}^N\left(\n a_{M,0,\ldots,0,0,0,\ldots,0}^{k,l}\n^2+\n a_{0,0,\ldots,0,M,0,\ldots,0}^{k,l}\n^2\right),
			\end{aligned}
		\end{equation}
		so $a_{M,0,\ldots,0,0,0,\ldots,0}^{k,l}=a_{0,0,\ldots,0,M,0,\ldots,0}^{k,l}=0$, $1\leq k,l\leq N$. In fact by considering $z_j^M\bz_j^M$, $j=1,2,\ldots,n$, we have $a_{0,\ldots,0,M,0,\ldots,0}^{k,l}=0$ where $M$ at the $j$-th and $(n+j)$-th position, $1\leq k,l\leq N$. Next consider the coefficient of $z_1^{M-1}z_2\bz_1^{M-1}\bz_2$, we have
		\begin{equation}
			\begin{aligned}
				0=&\sum_{l=1}^N\left(a_{M-1,1,\ldots,0,0,0,\ldots,0}^{k,l}\overline{a_{M-1,1,\ldots,0,0,0,\ldots,0}^{k,l}}\l_{2,1}^{M-1}+a_{0,0,\ldots,0,M-1,1,\ldots,0}^{k,l}\overline{a_{0,0,\ldots,0,M-1,1,\ldots,0}^{k,l}}\l_{2,1}^{M-1}\right)\\
				=&\l_{2,1}^{M-1}\sum_{l=1}^N\left(\n a_{M-1,1,\ldots,0,0,0,\ldots,0}^{k,l}\n^2+\n a_{0,0,\ldots,0,M-1,1,\ldots,0}^{k,l}\n^2\right),
			\end{aligned}
		\end{equation}
		so $a_{M-1,1,\ldots,0,0,0,\ldots,0}^{k,l}=a_{0,0,\ldots,0,M-1,1,\ldots,0}^{k,l}=0$, $1\leq k,l\leq N$. Similarly $a_{0,\ldots,M-1,\ldots,0,1,0,\ldots,0}^{k,l}=0$ where $M-1$ at the $s$-th and $(n+s)$-th position, 1 at the $t$-th and $(n+t)$-th position, $1\leq s,t\leq n$, $s\neq t$, $1\leq k,l\leq N$. Continue the above procedure, the formula \eqref{pqm} guarantees that we can prove $a_{p_1,\ldots,p_{n},q_1,\ldots,q_{n}}^{k,l}=0$ for all $(p_1,\ldots,p_{n},q_1,\ldots,q_{n})\in\N^{2n}$ satisfying $\sum_{j=1}^n(p_j+q_j)=M$ and $1\leq k,l\leq N$, but this contradicts the definition of $M$. So we must have $M=0$ and hence $p_{k,l}\in\C$ for all $1\leq k,l\leq N$, which implies $\P\in P(\C)$.
	\end{proof}
	For $\cron$, since $x$ commutes with everything, by replacing all $a_{p_1,\ldots,p_{n},q_1,\ldots,q_{n}}^{k,l}$ with polynomials of $x$ in the proof of Proposition \ref{p2n}, we have the following similar result.
	\begin{proposition}\label{p2n1}
		If $\P\in P(\cron)$, then $\P\in P(\C)$.
	\end{proposition}
	Then by the definition of $K_0$ groups, we have
	\begin{theorem}\label{kns}
		$K_0(\crm)=K_0(\C)=\Z$ for $\forall m\in\N^*$ and $\Theta$.
	\end{theorem}

	\section{The smooth cases}
	In this section we compute $K_0(\incrm)$. According to the parity of $m$, we divide it into two parts. Let's first consider the case where $m$ is even, $m=2n$. In the smooth cases, there's no longer result like Proposition \ref{p2n} and \ref{p2n1}. But by constructing a suitable unitary matrix $\mathcal{U}$, for any $\mathcal{P}\in P(\incren)$, we show that there is an $r\in\N$, depends only on $\mathcal{P}$, such that
	\begin{equation}
		\mathcal{U}\mathcal{P}\mathcal{U}^*=\begin{pmatrix}
			I_r & 0\\
			0 & 0\\
		\end{pmatrix},
	\end{equation}	
	which implies
 	\begin{theorem}\label{inp2n}
		$K_0(\incren)=\Z$ for $\forall n\in\N^*$ and $\Theta$.
	\end{theorem}
	\begin{proof}
		Like before, let $$\P=\begin{pmatrix}
			p_{1,1}&\ldots&p_{1,N}\\
			\vdots&\ddots&\vdots\\
			p_{N,1}&\ldots&p_{N,N}\\
		\end{pmatrix}\in M_N(C^\infty_{\mathrm{alg}}(\R_{\Theta}^{2n}))$$
		be an $N$-dimensional projector of $C^\infty_{\mathrm{alg}}(\R_{\Theta}^{2n})$. Then $\P^2=\P=\P^*$, and
		\begin{align*}
			\begin{pmatrix}
				p_{1,1}&\ldots&p_{1,N}\\
				\vdots&\ddots&\vdots\\
				p_{N,1}&\ldots&p_{N,N}\\
			\end{pmatrix}=\P=\P^2&=\P\P^*\\
			&=\begin{pmatrix}
				p_{1,1}&\ldots&p_{1,N}\\
				\vdots&\ddots&\vdots\\
				p_{N,1}&\ldots&p_{N,N}\\
			\end{pmatrix}\cdot\begin{pmatrix}
				p_{1,1}^*&\ldots&p_{N,1}^*\\
				\vdots&\ddots&\vdots\\
				p_{1,N}^*&\ldots&p_{N,N}^*\\
			\end{pmatrix}\\
			&=\begin{pmatrix}
				\sum_{k=1}^Np_{1,k}p_{1,k}^*&\ldots&\sum_{k=1}^Np_{1,k}p_{N,k}^*\\
				\vdots&\ddots&\vdots\\
				\sum_{k=1}^Np_{N,k}p_{1,k}^*&\ldots&\sum_{k=1}^Np_{N,k}p_{N,k}^*\\
			\end{pmatrix},
		\end{align*}
		so for $1\leq k,l\leq N$,
		\begin{equation}\label{pk}
			p_{k,l}=\sum_{j=1}^Np_{k,j}p_{l,j}^*.
		\end{equation}
		Let
		\begin{equation}
			p_{k,l}=\sum_{p_1,\ldots,p_{n},q_1,\ldots,q_{n}=0}^\infty a_{p_1,\ldots,p_{n},q_1,\ldots,q_{n}}^{k,l}z_1^{p_1}\ldots z_n^{p_n}\bz_{1}^{q_{1}}\ldots\bz_{n}^{q_{n}},~a_{p_1,\ldots,p_{n},q_1,\ldots,q_{n}}^{k,l}\in\C,~1\leq k,l\leq N,
		\end{equation}
		then 
		\begin{align*}
			&\sum_{j=1}^Np_{k,j}p_{l,j}^*\\
			=&\sum_{j=1}^N\sum_{p_1,\ldots,p_{n},q_1,\ldots,q_{n}=0}^\infty\sum_{p'_1,\ldots,p'_{n},q'_1,\ldots,q'_{n}=0}^\infty a_{p_1,\ldots,p_{n},q_1,\ldots,q_{n}}^{k,j}\overline{a_{p'_1,\ldots,p'_{n},q'_1,\ldots,q'_{n}}^{l,j}}z_1^{p_1}\ldots z_n^{p_n}\bz_{1}^{q_{1}}\ldots\bz_{n}^{q_{n}}z_n^{q'_n}\ldots z_1^{q'_1}\\
			&\times\bz_n^{p'_n}\ldots\bz_1^{p'_1}\\
			=&\sum_{j=1}^N\sum_{p_1,\ldots,p_{n},q_1,\ldots,q_{n}=0}^\infty\sum_{p'_1,\ldots,p'_{n},q'_1,\ldots,q'_{n}=0}^\infty a_{p_1,\ldots,p_{n},q_1,\ldots,q_{n}}^{k,j}\overline{a_{p'_1,\ldots,p'_{n},q'_1,\ldots,q'_{n}}^{l,j}}\left(\prod_{s=1}^{n-1}\prod_{r=s+1}^n\l_{r,s}^{p'_rp'_s+q'_rq'_s}\right)\\
			&\times z_1^{p_1}\ldots z_n^{p_n}\bz_{1}^{q_{1}}\ldots\bz_{n}^{q_{n}}z_1^{q'_1}\ldots z_n^{q'_n}\bz_{1}^{p'_{1}}\ldots\bz_{n}^{p'_{n}}\\
			=&\sum_{j=1}^N\sum_{p_1,\ldots,p_{n},q_1,\ldots,q_{n}=0}^\infty\sum_{p'_1,\ldots,p'_{n},q'_1,\ldots,q'_{n}=0}^\infty a_{p_1,\ldots,p_{n},q_1,\ldots,q_{n}}^{k,j}\overline{a_{p'_1,\ldots,p'_{n},q'_1,\ldots,q'_{n}}^{l,j}}\left(\prod_{s=1}^{n-1}\prod_{r=s+1}^n\l_{r,s}^{p'_rp'_s+q'_rq'_s}\right)\\
			&\times\prod_{l=1}^n\left(\prod_{k=l}^n\l_{r,s}^{p_rq'_s+q_rp'_s}\prod_{k=1}^n\l_{l,k}^{q_rq'_s}\right)z_1^{p_1+q'_1}\ldots z_n^{p_n+q'_n}\bz_{1}^{q_{1}+p'_1}\ldots\bz_{n}^{q_{n}+p'_n}\\
			=&\sum_{j=1}^N\sum_{p_1,\ldots,p_{n},q_1,\ldots,q_{n}=0}^\infty\sum_{p'_1,\ldots,p'_{n},q'_1,\ldots,q'_{n}=0}^\infty a_{p_1,\ldots,p_{n},q_1,\ldots,q_{n}}^{k,j}\overline{a_{p'_1,\ldots,p'_{n},q'_1,\ldots,q'_{n}}^{l,j}}z_1^{p_1+q'_1}\ldots z_n^{p_n+q'_n}\bz_{1}^{q_{1}+p'_1}\ldots\bz_{n}^{q_{n}+p'_n}\\
			&\times\prod_{s=1}^{n-1}\prod_{r=s+1}^n\l_{r,s}^{p'_rp'_s+q'_rq'_s+p_rq'_s+q_rp'_s+q_sq'_r-q_rq'_s},
		\end{align*}
		so we have
		\begin{equation}
			\begin{aligned}
				a_{p''_1,\ldots,p''_{n},q''_1,\ldots,q''_{n}}^{k,l}=&\sum_{j=1}^N\sum_{p_1+q'_1=p''_1,\ldots,p_{n}+q'_n=p''_n\atop q_1+p'_1=q''_1,\ldots,q_n+p'_n=q''_n}\left(\prod_{s=1}^{n-1}\prod_{r=s+1}^n\l_{r,s}^{p'_rp'_s+q'_rq'_s+p_rq'_s+q_rp'_s+q_sq'_r-q_rq'_s}\right)a_{p_1,\ldots,p_{n},q_1,\ldots,q_{n}}^{k,j}\\ 
				&\times \overline{a_{p'_1,\ldots,p'_{n},q'_1,\ldots,q'_{n}}^{l,j}}.
			\end{aligned}
		\end{equation}

		$\P=\P^*$ implies $p_{k,l}=p_{l,k}^*$, which means
		\begin{align*}
			&\sum_{p_1,\ldots,p_{n},q_1,\ldots,q_{n}=0}^\infty a_{p_1,\ldots,p_{n},q_1,\ldots,q_{n}}^{k,l}z_1^{p_1}\ldots z_n^{p_n}\bz_{1}^{q_{1}}\ldots\bz_{n}^{q_{n}}\\
			=&p_{k,l}=p_{l,k}^*\\
			=&\sum_{p_1,\ldots,p_{n},q_1,\ldots,q_{n}=0}^\infty \overline{a_{p_1,\ldots,p_{n},q_1,\ldots,q_{n}}^{l,k}}\left(z_1^{p_1}\ldots z_n^{p_n}\bz_{1}^{q_{1}}\ldots\bz_{n}^{q_{n}}\right)^*\\
			=&\sum_{p_1,\ldots,p_{n},q_1,\ldots,q_{n}=0}^\infty \overline{a_{p_1,\ldots,p_{n},q_1,\ldots,q_{n}}^{l,k}}z_n^{q_n}\ldots z_1^{q_1}\bz_n^{p_n}\ldots\bz_1^{p_1}\\
			=&\sum_{p_1,\ldots,p_{n},q_1,\ldots,q_{n}=0}^\infty \overline{a_{p_1,\ldots,p_{n},q_1,\ldots,q_{n}}^{l,k}}\left(\prod_{s=1}^{n-1}\prod_{r=s+1}^n\l_{r,s}^{p_rp_s+q_rq_s}\right)z_1^{q_1}\ldots z_n^{q_n}\bz_{1}^{p_{1}}\ldots\bz_{n}^{p_{n}},
		\end{align*}
		so we have
		\begin{equation}\label{self}
			a_{p_1,\ldots,p_{n},q_1,\ldots,q_{n}}^{k,l}=\overline{a_{q_1,q_2,\ldots,q_{n},p_1,p_2,\ldots,p_{n}}^{l,k}}\left(\prod_{s=1}^{n-1}\prod_{r=s+1}^n\l_{r,s}^{p_rp_s+q_rq_s}\right),
		\end{equation}
		and therefore  we also have
		\begin{equation}\label{aaa}
			\begin{aligned}
				a_{p''_1,\ldots,p''_{n},q''_1,\ldots,q''_{n}}^{k,l}=&\sum_{j=1}^N\sum_{p_1+p'_1=p''_1,\ldots,p_{n}+p'_n=p''_n\atop q_1+q'_1=q''_1,\ldots,q_n+q'_n=q''_n}\left(\prod_{s=1}^{n-1}\prod_{r=s+1}^n\l_{r,s}^{p_rp'_s+q_rq'_s+q_sp'_r-q_rp'_s}\right)a_{p_1,\ldots,p_{n},q_1,\ldots,q_{n}}^{k,j}\\ 
				&\times a_{p'_1,\ldots,p'_{n},q'_1,\ldots,q'_{n}}^{j,l}.
			\end{aligned}
		\end{equation}

		Note that for $1\leq k,l\leq N$,
		\begin{equation}
			a_{0,\ldots,0}^{k,l}=\overline{a_{0,\ldots,0}^{l,k}}
		\end{equation}
		and
		\begin{equation}
			a_{0,\ldots,0}^{k,l}=\sum_{j=1}^N a_{0,\ldots,0}^{k,j}a_{0,\ldots,0}^{j,l},
		\end{equation}
		so $\mathcal{A}=\{a_{0,\ldots,0}^{k,l}\}_{1\leq k,l\leq N}\in M_N(\C)$ is a projector. Denote $r:=\rank\left(\mathcal{A}\right),$
		the matrix rank of $\mathcal{A}$, then there is a unitary $N\times N$ matrix $\mathcal{Q}\in M_N(\C)$ such that
		\begin{equation}
			\mathcal{Q}\mathcal{B}\mathcal{Q}^*=\begin{pmatrix}
				I_r & 0\\
				0 & 0\\
			\end{pmatrix}.
		\end{equation}
		It's obvious that $\mathcal{Q}\mathcal{P}\mathcal{Q}^*$ is still a projector, so in fact we can assume
		\begin{equation}\label{asa}
			\{a_{0,\ldots,0}^{k,l}\}_{1\leq k,l\leq N}=\begin{pmatrix}
				I_r & 0\\
				0 & 0\\
			\end{pmatrix}.
		\end{equation}
		
		Next we show that there is a unitary matrix $\mathcal{U}\in M_N(C^\infty_{\mathrm{alg}}(\R_{\Theta}^{2n}))$ such that
		\begin{equation}
			\mathcal{U}\mathcal{P}\mathcal{U}^*=\begin{pmatrix}
				I_r & 0\\
				0 & 0\\
			\end{pmatrix}.
		\end{equation}	
		Let $$\mathcal{U}=\begin{pmatrix}
			u_{1,1}&\ldots&u_{1,N}\\
			\vdots&\ddots&\vdots\\
			u_{N,1}&\ldots&u_{N,N}\\
		\end{pmatrix}\in M_N(C^\infty_{\mathrm{alg}}(\R_{\Theta}^{2n}))$$
		be a unitary matrix where $u_{k,l}\in\incren$, $1\leq k,l\leq N$. Then $\mathcal{U}\mathcal{U}^*=I_N$ implies
		\begin{equation}
			\sum_{j=1}^Nu_{k,j}u_{l,j}^*= \delta_{k,l}
		\end{equation}
		where $\delta_{k,l}$ is the Kronecker delta. For $1\leq k,l\leq N$, let
		\begin{equation}
			u_{k,l}=\sum_{p_1,\ldots,p_{n},q_1,\ldots,q_{n}=0}^\infty v_{p_1,\ldots,p_{n},q_1,\ldots,q_{n}}^{k,l}z_1^{p_1}\ldots z_n^{p_n}\bz_{1}^{q_{1}}\ldots\bz_{n}^{q_{n}},
		\end{equation}
		and for all $(p_1,\ldots,p_{n},q_1,\ldots,q_{n})\in\N^{2n}$, for convenience we  define
		\begin{equation}
			w_{p_1,\ldots,p_{n},q_1,\ldots,q_{n}}^{k,l}:=\overline{v_{q_1,\ldots,q_{n},p_1,\ldots,p_{n}}^{l,k}},
		\end{equation}
		then 
		\begin{align*}
			&\sum_{j=1}^Nu_{k,j}u_{l,j}^*\\
			=&\sum_{j=1}^N\sum_{p_1,\ldots,p_{n},q_1,\ldots,q_{n}=0}^\infty\sum_{p'_1,\ldots,p'_{n},q'_1,\ldots,q'_{n}=0}^\infty v_{p_1,\ldots,p_{n},q_1,\ldots,q_{n}}^{k,j}\overline{v_{p'_1,\ldots,p'_{n},q'_1,\ldots,q'_{n}}^{l,j}}z_1^{p_1}\ldots z_n^{p_n}\bz_{1}^{q_{1}}\ldots\bz_{n}^{q_{n}}z_n^{q'_n}\ldots z_1^{q'_1}\\
			&\times\bz_n^{p'_n}\ldots\bz_1^{p'_1}\\
			=&\sum_{j=1}^N\sum_{p_1,\ldots,p_{n},q_1,\ldots,q_{n}=0}^\infty\sum_{p'_1,\ldots,p'_{n},q'_1,\ldots,q'_{n}=0}^\infty v_{p_1,\ldots,p_{n},q_1,\ldots,q_{n}}^{k,j}\overline{v_{p'_1,\ldots,p'_{n},q'_1,\ldots,q'_{n}}^{l,j}}z_1^{p_1+q'_1}\ldots z_n^{p_n+q'_n}\bz_{1}^{q_{1}+p'_1}\ldots\bz_{n}^{q_{n}+p'_n}\\
			&\times\prod_{s=1}^{n-1}\prod_{r=s+1}^n\l_{r,s}^{p'_rp'_s+q'_rq'_s+p_rq'_s+q_rp'_s+q_sq'_r-q_rq'_s}\\
			=&\sum_{j=1}^N\sum_{p_1,\ldots,p_{n},q_1,\ldots,q_{n}=0}^\infty\sum_{p'_1,\ldots,p'_{n},q'_1,\ldots,q'_{n}=0}^\infty v_{p_1,\ldots,p_{n},q_1,\ldots,q_{n}}^{k,j}w_{p'_1,\ldots,p'_{n},q'_1,\ldots,q'_{n}}^{j,l}z_1^{p_1+p'_1}\ldots z_n^{p_n+p'_n}\bz_{1}^{q_{1}+q'_1}\ldots\bz_{n}^{q_{n}+q'_n}\\
			&\times\prod_{s=1}^{n-1}\prod_{r=s+1}^n\l_{r,s}^{p'_rp'_s+q'_rq'_s+p_rp'_s+q_rq'_s+q_sp'_r-q_rp'_s},
		\end{align*}
		so for $(p''_1,\ldots,p''_{n},q''_1,\ldots,q''_{n})\in\N^{2n}$, we have
		\begin{equation}\label{vw0}
			\begin{aligned}
				&\sum_{j=1}^N\sum_{p_1+p'_1=p''_1,\ldots,p_{n}+p'_n=p''_n\atop q_1+q'_1=q''_1,\ldots,q_n+q'_n=q''_n}v_{p_1,\ldots,p_{n},q_1,\ldots,q_{n}}^{k,j}w_{p'_1,\ldots,p'_{n},q'_1,\ldots,q'_{n}}^{j,l}\prod_{s=1}^{n-1}\prod_{r=s+1}^n\l_{r,s}^{p'_rp'_s+q'_rq'_s+p_rp'_s+q_rq'_s+q_sp'_r-q_rp'_s}\\
				=&\delta_{k,l}\delta_{\sum_{j=1}^N\left(p''_j+q''_j\right),0}.
			\end{aligned}
		\end{equation}

		Note that
		\begin{equation}
			\sum_{j=1}^N v_{0,\ldots,0}^{k,j}\overline{v_{0,\ldots,0}^{l,j}}=\delta_{k,l},
		\end{equation}
		so $\mathcal{V}:=\{v_{0,\ldots,0}^{k,l}\}_{1\leq k,l\leq N}\in M_N(\C)$ is a unitary $N\times N$ matrix. It's obvious that $\mathcal{U}\mathcal{V}^*$ is still unitary, so in fact we can assume 
		$$v_{0,\ldots,0}^{k,l}=\delta_{k,l}.$$
		Then the element in the $(k,l)$-position of $\mathcal{U}\mathcal{P}\mathcal{U}^*$ is		
		\begin{align*}
			&q_{k,l}\\
			:=&\sum_{p'''_1,\ldots,p'''_{n},q'''_1,\ldots,q'''_{n}=0}^\infty c_{p'''_1,\ldots,p'''_{n},q'''_1,\ldots,q'''_{n}}^{k,l}z_1^{p'''_1}\ldots z_n^{p'''_n}\bz_{1}^{q'''_{1}}\ldots\bz_{n}^{q'''_{n}}\\
			:=&\sum_{\a,\b=1}^Nu_{k,\a}p_{\a,\b}u_{l,\b}^*\\
			=&\sum_{\a,\b=1}^N\sum_{p_1,\ldots,p_{n},q_1,\ldots,q_{n}=0}^\infty v_{p_1,\ldots,p_{n},q_1,\ldots,q_{n}}^{k,\a}z_1^{p_1}\ldots z_n^{p_n}\bz_{1}^{q_{1}}\ldots\bz_{n}^{q_{n}}\sum_{p'_1,\ldots,p'_{n},q'_1,\ldots,q'_{n}=0}^\infty a_{p'_1,\ldots,p'_{n},q'_1,\ldots,q'_{n}}^{\a,\b}\\
			&\times z_1^{p'_1}\ldots z_n^{p'_n}\bz_{1}^{q'_{1}}\ldots\bz_{n}^{q'_{n}}\sum_{p''_1,\ldots,p''_{n},q''_1,\ldots,q''_{n}=0}^\infty \overline{v_{p''_1,\ldots,p''_{n},q''_1,\ldots,q''_{n}}^{l,\b}}z_n^{q''_n}\ldots z_1^{q''_1}\bz_n^{p''_n}\ldots\bz_1^{p''_1}\\
			=&\sum_{\a,\b=1}^N\sum_{{p_1,\ldots,p_{n},q_1,\ldots,q_{n}=0\atop p'_1,\ldots,p'_{n},q'_1,\ldots,q'_{n}=0}\atop p''_1,\ldots,p''_{n},q''_1,\ldots,q''_{n}=0}^\infty v_{p_1,\ldots,p_{n},q_1,\ldots,q_{n}}^{k,\a}a_{p'_1,\ldots,p'_{n},q'_1,\ldots,q'_{n}}^{\a,\b}\overline{v_{p''_1,\ldots,p''_{n},q''_1,\ldots,q''_{n}}^{l,\b}}\\
			&\times\left(\prod_{s=1}^{n-1}\prod_{r=s+1}^n\l_{r,s}^{p_rp'_s+q_rq'_s+q_sp'_r-q_rp'_s}\right)\left(\prod_{s=1}^{n-1}\prod_{r=s+1}^n\l_{r,s}^{p''_rp''_s+q''_rq''_s}\right)z_1^{p_1+p'_1}\ldots z_n^{p_n+p'_n}\bz_1^{q_1+q'_1}\ldots \bz_n^{q_n+q'_n}\\
			&\times z_1^{q''_1}\ldots z_n^{q''_n}\bz_{1}^{p''_{1}}\ldots\bz_{n}^{p''_{n}}\\
			=&\sum_{\a,\b=1}^N\sum_{{p_1,\ldots,p_{n},q_1,\ldots,q_{n}=0\atop p'_1,\ldots,p'_{n},q'_1,\ldots,q'_{n}=0}\atop p''_1,\ldots,p''_{n},q''_1,\ldots,q''_{n}=0}^\infty v_{p_1,\ldots,p_{n},q_1,\ldots,q_{n}}^{k,\a}a_{p'_1,\ldots,p'_{n},q'_1,\ldots,q'_{n}}^{\a,\b}\overline{v_{p''_1,\ldots,p''_{n},q''_1,\ldots,q''_{n}}^{l,\b}}\\
			&\times\left(\prod_{s=1}^{n-1}\prod_{r=s+1}^n\l_{r,s}^{p_rp'_s+q_rq'_s+q_sp'_r-q_rp'_s}\right)\left(\prod_{s=1}^{n-1}\prod_{r=s+1}^n\l_{r,s}^{p''_rp''_s+q''_rq''_s}\right)\\
			&\times\left(\prod_{s=1}^{n-1}\prod_{r=s+1}^n\l_{r,s}^{\left(p_r+p'_r\right)q''_s+\left(q_r+q'_r\right)p''_s+\left(q_s+q'_s\right)q''_r-\left(q_r+q'_r\right)q''_s}\right)z_1^{p_1+p'_1+q''_1}\ldots z_n^{p_n+p'_n+q''_n}\\
			&\times\bz_1^{q_1+q'_1+p''_1}\ldots \bz_n^{q_n+q'_n+p''_n}\\
			=&\sum_{\a,\b=1}^N\sum_{{p_1,\ldots,p_{n},q_1,\ldots,q_{n}=0\atop p'_1,\ldots,p'_{n},q'_1,\ldots,q'_{n}=0}\atop p''_1,\ldots,p''_{n},q''_1,\ldots,q''_{n}=0}^\infty v_{p_1,\ldots,p_{n},q_1,\ldots,q_{n}}^{k,\a}a_{p'_1,\ldots,p'_{n},q'_1,\ldots,q'_{n}}^{\a,\b}\overline{v_{p''_1,\ldots,p''_{n},q''_1,\ldots,q''_{n}}^{l,\b}}\\
			&\times\prod_{s=1}^{n-1}\prod_{r=s+1}^n\l_{r,s}^{p_rp'_s+q_rq'_s+q_sp'_r-q_rp'_s+\left(p_r+p'_r+q''_r\right)q''_s+\left(q_r+q'_r+p''_r\right)p''_s+\left(q_s+q'_s\right)q''_r-\left(q_r+q'_r\right)q''_s}\\
			&\times z_1^{p_1+p'_1+q''_1}\ldots z_n^{p_n+p'_n+q''_n}\bz_1^{q_1+q'_1+p''_1}\ldots \bz_n^{q_n+q'_n+p''_n}\\
			=&\sum_{\a,\b=1}^N\sum_{p_1+p'_1+p''_1=p'''_1,\ldots,p_n+p'_n+p''_n=p'''_n\atop q_1+q'_1+q''_1=q'''_1,\ldots,q_n+q'_n+q''_n=q'''_n}^\infty v_{p_1,\ldots,p_{n},q_1,\ldots,q_{n}}^{k,\a}a_{p'_1,\ldots,p'_{n},q'_1,\ldots,q'_{n}}^{\a,\b}w_{p''_1,\ldots,p''_{n},q''_1,\ldots,q''_{n}}^{\b,l}\\
			&\times\prod_{s=1}^{n-1}\prod_{r=s+1}^n\l_{r,s}^{p_rp'_s+q_rq'_s+q_sp'_r-q_rp'_s+\left(p_r+p'_r+p''_r\right)p''_s+\left(q_r+q'_r+q''_r\right)q''_s+\left(q_s+q'_s\right)p''_r-\left(q_r+q'_r\right)p''_s}\\
			&\times z_1^{p_1+p'_1+p''_1}\ldots z_n^{p_n+p'_n+p''_n}\bz_1^{q_1+q'_1+q''_1}\ldots \bz_n^{q_n+q'_n+q''_n}\\
			=&\sum_{\a,\b=1}^N\sum_{p_1+p'_1+p''_1=p'''_1,\ldots,p_n+p'_n+p''_n=p'''_n\atop q_1+q'_1+q''_1=q'''_1,\ldots,q_n+q'_n+q''_n=q'''_n}^\infty v_{p_1,\ldots,p_{n},q_1,\ldots,q_{n}}^{k,\a}a_{p'_1,\ldots,p'_{n},q'_1,\ldots,q'_{n}}^{\a,\b}w_{p''_1,\ldots,p''_{n},q''_1,\ldots,q''_{n}}^{\b,l}\\
			&\times\prod_{s=1}^{n-1}\prod_{r=s+1}^n\l_{r,s}^{p_rp'_s+q_rq'_s+q_sp'_r-q_rp'_s+p'''_rp''_s+q'''_rq''_s+\left(q'''_s-q''_s\right)p''_r-\left(q'''_r-q''_r\right)p''_s}\\
			&\times z_1^{p_1+p'_1+p''_1}\ldots z_n^{p_n+p'_n+p''_n}\bz_1^{q_1+q'_1+q''_1}\ldots \bz_n^{q_n+q'_n+q''_n}.
		\end{align*}
		The coefficient of $z_1$ in $q_{k,l}$ is
		\begin{align*}
			c_{1,\ldots,0}^{k,l}=&\sum_{\a,\b=1}^N\left(v_{1,\ldots,0}^{k,\a}a_{0,\ldots,0}^{\a,\b}w_{0,\ldots,0}^{\b,l}+v_{0,\ldots,0}^{k,\a}a_{1,\ldots,0}^{\a,\b}w_{0,\ldots,0}^{\b,l}+v_{0,\ldots,0}^{k,\a}a_{0,\ldots,0}^{\a,\b}w_{1,\ldots,0}^{\b,l}\right)\\
			=&v_{1,\ldots,0}^{k,l}a_{0,\ldots,0}^{l,l}+a_{1,\ldots,0}^{k,l}+a_{0,\ldots,0}^{k,k}w_{1,\ldots,0}^{k,l}.
		\end{align*}
		Note that
		\begin{equation}
			a_{0,\ldots,1,\ldots,0}^{k,l}=\sum_{j=1}^N\left(a_{0,\ldots,1,\ldots,0}^{k,j}a_{0,\ldots,0}^{j,l}+a_{0,\ldots,0}^{k,j}a_{0,\ldots,1,\ldots,0}^{j,l}\right)=a_{0,\ldots,1,\ldots,0}^{k,l}\left(a_{0,\ldots,0}^{k,k}+a_{0,\ldots,0}^{l,l}\right),
		\end{equation}
		where 1 is at the $h$-th position, so $a_{0,\ldots,1,\ldots,0}^{k,l}=0$ where $1\leq k,l\leq r$ and $r+1\leq k,l\leq N$, $1\leq h\leq 2n$. Also,
		\begin{equation}
			v_{1,\ldots,0}^{k,l}+w_{1,\ldots,0}^{k,l}=0
		\end{equation}
		where 1 is at the $(n+1)$-th position. If we set
		\begin{equation}
			v_{1,\ldots,0}^{k,l}=\left\{\begin{aligned}
				a_{1,\ldots,0}^{k,l},~&1\leq k\leq r<l\leq N,\\
				-a_{1,\ldots,0}^{k,l},~&1\leq l\leq r<k\leq N,
			\end{aligned}\right.
		\end{equation}
		then $c_{1,\ldots,0}^{k,l}=0$ for all  $1\leq k,l\leq N$, which means the coefficient of $z_1$ in $q_{k,l}$ is 0. In fact this implies the coefficient of $\bz_1$ in $q_{k,l}$ is also 0. Similarly for $1\leq h\leq n$, set
		\begin{equation}
			v_{0,\ldots,1,\ldots,0}^{k,l}=\left\{\begin{aligned}
				a_{0,\ldots,1,\ldots,0}^{k,l},~&1\leq k\leq r<l\leq N,\\
				-a_{0,\ldots,1,\ldots,0}^{k,l},~&1\leq l\leq r<k\leq N
			\end{aligned}\right.
		\end{equation}
		where all 1's are at the $h$-th position, the coefficient of $z_h$ and $\bz_h$ in $q_{k,l}$ are 0 for all $1\leq k,l\leq N$.

		Next,
		\begin{align*}
			a_{0,\ldots,2,\ldots,0}^{k,l}=&\sum_{j=1}^N\left(a_{0,\ldots,2,\ldots,0}^{k,j}a_{0,\ldots,0}^{j,l}+a_{0,\ldots,1,\ldots,0}^{k,j}a_{0,\ldots,1,\ldots,0}^{j,l}+a_{0,\ldots,0}^{k,j}a_{0,\ldots,2,\ldots,0}^{j,l}\right)\\
			=&\left\{\begin{aligned}
				2a_{0,\ldots,2,\ldots,0}^{k,l}+\sum_{j=r+1}^Na_{0,\ldots,1,\ldots,0}^{k,j}a_{0,\ldots,1,\ldots,0}^{j,l},~&1\leq k,l\leq r,\\
				\sum_{j=1}^r a_{0,\ldots,1,\ldots,0}^{k,j}a_{0,\ldots,1,\ldots,0}^{j,l},~&r+1\leq k,l\leq N,\\
				a_{0,\ldots,2,\ldots,0}^{k,l},&\text{ other cases,}
			\end{aligned}\right.
		\end{align*}
		where 1 and 2 are at the $h$-th position,   so
		\begin{equation}
			a_{0,\ldots,2,\ldots,0}^{k,l}=\left\{\begin{aligned}
				-\sum_{j=r+1}^Na_{0,\ldots,1,\ldots,0}^{k,j}a_{0,\ldots,1,\ldots,0}^{j,l},~&1\leq k,l\leq r,\\
				\sum_{j=1}^r a_{0,\ldots,1,\ldots,0}^{k,j}a_{0,\ldots,1,\ldots,0}^{j,l},~&r+1\leq k,l\leq N.
			\end{aligned}\right.
		\end{equation}
		Also,
		\begin{equation}
			v_{2,\ldots,0}^{k,l}+\sum_{j=1}^Nv_{1,\ldots,0}^{k,j}w_{1,\ldots,0}^{j,l}+w_{2,\ldots,0}^{k,l}=0.
		\end{equation}
		The coefficient of $z_1^2$ in $q_{k,l}$ is
		\begin{align*}
			c_{2,\ldots,0}^{k,l}=&\sum_{\a,\b=1}^N\left(v_{2,\ldots,0}^{k,\a}a_{0,\ldots,0}^{\a,\b}w_{0,\ldots,0}^{\b,l}+v_{0,\ldots,0}^{k,\a}a_{2,\ldots,0}^{\a,\b}w_{0,\ldots,0}^{\b,l}+v_{0,\ldots,0}^{k,\a}a_{0,\ldots,0}^{\a,\b}w_{2,\ldots,0}^{\b,l}+v_{1,\ldots,0}^{k,\a}a_{1,\ldots,0}^{\a,\b}w_{0,\ldots,0}^{\b,l}\right.\\
			&\left.\quad+v_{1,\ldots,0}^{k,\a}a_{0,\ldots,0}^{\a,\b}w_{1,\ldots,0}^{\b,l}+v_{0,\ldots,0}^{k,\a}a_{1,\ldots,0}^{\a,\b}w_{1,\ldots,0}^{\b,l}\right)\\
			=&v_{2,\ldots,0}^{k,l}a_{0,\ldots,0}^{l,l}+a_{2,\ldots,0}^{k,l}+a_{0,\ldots,0}^{k,k}w_{2,\ldots,0}^{k,l}+\sum_{\a=1}^Nv_{1,\ldots,0}^{k,\a}a_{1,\ldots,0}^{\a,l}+\sum_{\a=1}^Nv_{1,\ldots,0}^{k,\a}a_{0,\ldots,0}^{\a,\a}w_{1,\ldots,0}^{\a,l}+\sum_{\a=1}^Na_{1,\ldots,0}^{k,\a}w_{1,\ldots,0}^{\a,l}\\
			=&\left\{\begin{aligned}
				0,~&1\leq k,l\leq r,\\
				0,~&r+1\leq k,l\leq N,\\
				a_{2,\ldots,0}^{k,l}-v_{2,\ldots,0}^{k,l}+\sum_{j=1}^rv_{1,\ldots,0}^{k,j}a_{1,\ldots,0}^{j,l},~&1\leq k\leq r<l\leq N,\\
				a_{2,\ldots,0}^{k,l}+v_{2,\ldots,0}^{k,l}+\sum_{j=r+1}^Nv_{1,\ldots,0}^{k,j}a_{1,\ldots,0}^{j,l},~&1\leq l\leq r<k\leq N.
			\end{aligned}\right.
		\end{align*}
		If we set
		\begin{equation}
			v_{2,\ldots,0}^{k,l}=\left\{\begin{aligned}
				a_{2,\ldots,0}^{k,l}+\sum_{j=1}^rv_{1,\ldots,0}^{k,j}a_{1,\ldots,0}^{j,l},~&1\leq k\leq r<l\leq N,\\
				-a_{2,\ldots,0}^{k,l}-\sum_{j=r+1}^Nv_{1,\ldots,0}^{k,j}a_{1,\ldots,0}^{j,l},~&1\leq l\leq r<k\leq N,
			\end{aligned}\right.
		\end{equation}
		then $c_{2,\ldots,0}^{k,l}=0$ for all  $1\leq k,l\leq N$, which means the coefficient of $z_1^2$ and $\bz_1^2$ in $q_{k,l}$ are  0.  Like before, for $1\leq h\leq n$, set
		\begin{equation}
			v_{0,\ldots,2,\ldots,0}^{k,l}=\left\{\begin{aligned}
				a_{0,\ldots,2,\ldots,0}^{k,l}+\sum_{j=1}^rv_{0,\ldots,1,\ldots,0}^{k,j}a_{0,\ldots,1,\ldots,0}^{j,l},~&1\leq k\leq r<l\leq N,\\
				-a_{0,\ldots,2,\ldots,0}^{k,l}-\sum_{j=r+1}^Nv_{0,\ldots,1,\ldots,0}^{k,j}a_{0,\ldots,1,\ldots,0}^{j,l},~&1\leq l\leq r<k\leq N,
			\end{aligned}\right.
		\end{equation}
		where all 1's and 2's are at the $h$-th position, the coefficient of $z_h^2$ and $\bz_h^2$ in $q_{k,l}$ are 0 for all $1\leq k,l\leq N$. 
		
		In general, for all 
		$(p''_1,\ldots,p''_{n},q''_1,\ldots,q''_{n})\in\N^{2n}$ with $\sum_{j=1}^n\left(p''_j+q''_j\right)\geq 1$, we set
		\begin{equation}\label{vkrl}
			v_{p''_1,\ldots,p''_{n},q''_1,\ldots,q''_{n}}^{k,l}=\left\{\begin{aligned}
			a_{p''_1,\ldots,p''_{n},q''_1,\ldots,q''_{n}}^{k,l}+\sum_{j=1}^N\sum_{{p_1+p'_1=p''_1,\ldots,p_n+p'_n=p''_n\atop q_1+q'_1=q''_1,\ldots,q_n+q'_n=q''_n}\atop\sum_{j=1}^n\left(p_j+q_j\right),\sum_{j=1}^n\left(p'_j+q'_j\right)\geq 1}v_{p_1,\ldots,p_{n},q_1,\ldots,q_{n}}^{k,j}a_{p'_1,\ldots,p'_{n},q'_1,\ldots,q'_{n}}^{j,l}\\
			\times\prod_{s=1}^{n-1}\prod_{r=s+1}^n\l_{r,s}^{p_rp'_s+q_rq'_s+q_sp'_r-q_rp'_s},~1\leq k\leq r<l\leq N,\\
			-a_{p''_1,\ldots,p''_{n},q''_1,\ldots,q''_{n}}^{k,l}-\sum_{j=1}^N\sum_{{p_1+p'_1=p''_1,\ldots,p_n+p'_n=p''_n\atop q_1+q'_1=q''_1,\ldots,q_n+q'_n=q''_n}\atop\sum_{j=1}^n\left(p_j+q_j\right),\sum_{j=1}^n\left(p'_j+q'_j\right)\geq 1}v_{p_1,\ldots,p_{n},q_1,\ldots,q_{n}}^{k,j}a_{p'_1,\ldots,p'_{n},q'_1,\ldots,q'_{n}}^{j,l}\\
			\times\prod_{s=1}^{n-1}\prod_{r=s+1}^n\l_{r,s}^{p_rp'_s+q_rq'_s+q_sp'_r-q_rp'_s},~1\leq l\leq r<k\leq N,\\
			\end{aligned}\right.
		\end{equation}
		and for $1\leq k,l\leq r$ and $r+1\leq k,l\leq N$, $v_{p''_1,\ldots,p''_{n},q''_1,\ldots,q''_{n}}^{k,l}$ are chosen such that \eqref{vw0} are satisfied for all $(p''_1,\ldots,p''_{n},q''_1,\ldots,q''_{n})\in\N^{2n}$. It's not hard to see that this can be done from $\sum_{j=1}^n\left(p''_j+q''_j\right)=1$ to any far.

		For $1\leq k,l\leq r$, from \eqref{aaa} we have
		\begin{equation}\label{aklr}
			a_{p''_1,\ldots,p''_{n},q''_1,\ldots,q''_{n}}^{k,l}=-\sum_{{p_1+p'_1=p''_1,\ldots,p_n+p'_n=p''_n\atop q_1+q'_1=q''_1,\ldots,q_n+q'_n=q''_n }\atop\sum_{j=1}^n\left(p_j+q_j\right),\sum_{j=1}^n\left(p'_j+q'_j\right)\geq 1}a_{p_1,\ldots,p_{n},q_1,\ldots,q_{n}}^{k,j}a_{p'_1,\ldots,p'_{n},q'_1,\ldots,q'_{n}}^{j,l}\prod_{s=1}^{n-1}\prod_{r=s+1}^n\l_{r,s}^{p_rp'_s+q_rq'_s+q_sp'_r-q_rp'_s}.
		\end{equation}
	By applying \eqref{vkrl} and \eqref{aklr} repeatly, the sum

	\begin{align*}
		&\sum_{j=1}^N\sum_{p_1+p'_1=p''_1,\ldots,p_n+p'_n=p''_n\atop q_1+q'_1=q''_1,\ldots,q_n+q'_n=q''_n}v_{p_1,\ldots,p_{n},q_1,\ldots,q_{n}}^{k,j}a_{p'_1,\ldots,p'_{n},q'_1,\ldots,q'_{n}}^{j,l}\prod_{s=1}^{n-1}\prod_{r=s+1}^n\l_{r,s}^{p_rp'_s+q_rq'_s+q_sp'_r-q_rp'_s}\\
		=&\sum_{j=1}^N\sum_{{p_1+p'_1=p''_1,\ldots,p_n+p'_n=p''_n\atop q_1+q'_1=q''_1,\ldots,q_n+q'_n=q''_n }\atop\sum_{j=1}^n\left(p_j+q_j\right),\sum_{j=1}^n\left(p'_j+q'_j\right)\geq 1}v_{p_1,\ldots,p_{n},q_1,\ldots,q_{n}}^{k,j}a_{p'_1,\ldots,p'_{n},q'_1,\ldots,q'_{n}}^{j,l}\prod_{s=1}^{n-1}\prod_{r=s+1}^n\l_{r,s}^{p_rp'_s+q_rq'_s+q_sp'_r-q_rp'_s}\\
		&+v_{p''_1,\ldots,p''_{n},q''_1,\ldots,q''_{n}}^{k,l}+a_{p''_1,\ldots,p''_{n},q''_1,\ldots,q''_{n}}^{k,l}\\
		=&\sum_{j=1}^r\sum_{{p_1+p'_1=p''_1,\ldots,p_n+p'_n=p''_n\atop q_1+q'_1=q''_1,\ldots,q_n+q'_n=q''_n }\atop\sum_{j=1}^n\left(p_j+q_j\right),\sum_{j=1}^n\left(p'_j+q'_j\right)\geq 1}v_{p_1,\ldots,p_{n},q_1,\ldots,q_{n}}^{k,j}a_{p'_1,\ldots,p'_{n},q'_1,\ldots,q'_{n}}^{j,l}\prod_{s=1}^{n-1}\prod_{r=s+1}^n\l_{r,s}^{p_rp'_s+q_rq'_s+q_sp'_r-q_rp'_s}\\
		&+\sum_{j=r+1}^N\sum_{{p_1+p'_1=p''_1,\ldots,p_n+p'_n=p''_n\atop q_1+q'_1=q''_1,\ldots,q_n+q'_n=q''_n }\atop\sum_{j=1}^n\left(p_j+q_j\right),\sum_{j=1}^n\left(p'_j+q'_j\right)\geq 1}v_{p_1,\ldots,p_{n},q_1,\ldots,q_{n}}^{k,j}a_{p'_1,\ldots,p'_{n},q'_1,\ldots,q'_{n}}^{j,l}\prod_{s=1}^{n-1}\prod_{r=s+1}^n\l_{r,s}^{p_rp'_s+q_rq'_s+q_sp'_r-q_rp'_s}\\
		&+v_{p''_1,\ldots,p''_{n},q''_1,\ldots,q''_{n}}^{k,l}+a_{p''_1,\ldots,p''_{n},q''_1,\ldots,q''_{n}}^{k,l}
	\end{align*}
	can be finally expressed by the sum of type VA:
	\begin{equation}\label{vaa}
		\begin{aligned}
			&v_{p_{1,0},\ldots,p_{n,0},q_{1,0},\ldots,q_{n,0}}^{k,\a_1}\prod_{j=1}^{n_1}a_{p_{1,2j-1},\ldots,p_{n,2j-1},q_{1,2j-1},\ldots,q_{n,2j-1}}^{\a_j,\b_j}a_{p_{1,2j},\ldots,p_{n,2j},q_{1,2j},\ldots,q_{n,2j}}^{\b_j,\a_{j+1}}\Phi_{\{p_{r,s},q_{r,s}\}_{1\leq r\leq n,0\leq s\leq 2n_1}}
		\end{aligned}
	\end{equation}
	and type AA:
	\begin{equation}\label{aa}
		\begin{aligned}
			&\prod_{j=1}^{n_2}a_{p'_{1,2j-1},\ldots,p'_{n,2j-1},q'_{1,2j-1},\ldots,q'_{n,2j-1}}^{\g_j,\d_j}a_{p'_{1,2j},\ldots,p'_{n,2j},q'_{1,2j},\ldots,q'_{n,2j}}^{\d_j,\g_{j+1}}\Phi_{\{p'_{r,s},q'_{r,s}\}_{1\leq r\leq n,1\leq s\leq 2n_2}},
		\end{aligned}
	\end{equation}	
	here $n_1,n_2\in\N$, $\sum_{v=0}^{2n_1}\sum_{u=1}^n\left(p_{u,v}+q_{u,v}\right)=\sum_{v=1}^{2n_2}\sum_{u=1}^n\left(p'_{u,v}+q'_{u,v}\right)=\sum_{u=1}^n\left(p''_u+q''_u\right)$, $1\leq\a_u\leq r<\b_u\leq N$, $1\leq\g_v\leq r<\d_v\leq N$, $1\leq u\leq n_1$, $1\leq v\leq n_2$, $\g_1=k$, $\a_{n_1+1}=\g_{n_2+1}=l$, $\Phi_{\{p_{r,s},q_{r,s}\}_{1\leq r\leq n,0\leq s\leq 2n_j}}$ is a complex number with norm 1 which is determined by $p_{r,s}$ and $q_{r,s}$ where $1\leq r\leq n$, $0\leq s\leq 2n_j$, $j=1,2$. 
	
	From \eqref{aklr} the expression of
	\begin{equation}
		\sum_{j=1}^r\sum_{{p_1+p'_1=p''_1,\ldots,p_n+p'_n=p''_n\atop q_1+q'_1=q''_1,\ldots,q_n+q'_n=q''_n }\atop\sum_{j=1}^n\left(p_j+q_j\right),\sum_{j=1}^n\left(p'_j+q'_j\right)\geq 1}v_{p_1,\ldots,p_{n},q_1,\ldots,q_{n}}^{k,j}a_{p'_1,\ldots,p'_{n},q'_1,\ldots,q'_{n}}^{j,l}\prod_{s=1}^{n-1}\prod_{r=s+1}^n\l_{r,s}^{p_rp'_s+q_rq'_s+q_sp'_r-q_rp'_s}
	\end{equation}
	contains only terms of type VA and all are equiped with negative signs; From \eqref{vkrl} the expression of
	\begin{equation}
		\sum_{j=r+1}^N\sum_{{p_1+p'_1=p''_1,\ldots,p_n+p'_n=p''_n\atop q_1+q'_1=q''_1,\ldots,q_n+q'_n=q''_n }\atop\sum_{j=1}^n\left(p_j+q_j\right),\sum_{j=1}^n\left(p'_j+q'_j\right)\geq 1}v_{p_1,\ldots,p_{n},q_1,\ldots,q_{n}}^{k,j}a_{p'_1,\ldots,p'_{n},q'_1,\ldots,q'_{n}}^{j,l}\prod_{s=1}^{n-1}\prod_{r=s+1}^n\l_{r,s}^{p_rp'_s+q_rq'_s+q_sp'_r-q_rp'_s}
	\end{equation}
	contains both types and all are equiped with positive signs; Finally, from
	\eqref{aklr} again we can see that the expression of $a_{p''_1,\ldots,p''_{n},q''_1,\ldots,q''_{n}}^{k,l}$ contains only terms of type AA and all are equiped with negative signs. Hence we have
	\begin{equation}
		\sum_{j=1}^N\sum_{p_1+p'_1=p''_1,\ldots,p_n+p'_n=p''_n\atop q_1+q'_1=q''_1,\ldots,q_n+q'_n=q''_n}v_{p_1,\ldots,p_{n},q_1,\ldots,q_{n}}^{k,j}a_{p'_1,\ldots,p'_{n},q'_1,\ldots,q'_{n}}^{j,l}\prod_{s=1}^{n-1}\prod_{r=s+1}^n\l_{r,s}^{p_rp'_s+q_rq'_s+q_sp'_r-q_rp'_s}=v_{p''_1,\ldots,p''_{n},q''_1,\ldots,q''_{n}}^{k,l}.
	\end{equation}
	Similarly we can prove that for $r+1\leq k,l\leq N$,
	\begin{equation}
		\sum_{j=1}^N\sum_{p_1+p'_1=p''_1,\ldots,p_n+p'_n=p''_n\atop q_1+q'_1=q''_1,\ldots,q_n+q'_n=q''_n}v_{p_1,\ldots,p_{n},q_1,\ldots,q_{n}}^{k,j}a_{p'_1,\ldots,p'_{n},q'_1,\ldots,q'_{n}}^{j,l}\prod_{s=1}^{n-1}\prod_{r=s+1}^n\l_{r,s}^{p_rp'_s+q_rq'_s+q_sp'_r-q_rp'_s}=0.
	\end{equation}
In summary, for all  $(p'''_1,\ldots,p'''_{n},q'''_1,\ldots,q'''_{n})\in\N^{2n}$ with $\sum_{j=1}^n\left(p'''_j+q'''_j\right)\geq1$ we have
\begin{equation}
	\begin{aligned}
		&\sum_{j=1}^N\sum_{p_1+p'_1=p''_1,\ldots,p_n+p'_n=p''_n\atop q_1+q'_1=q''_1,\ldots,q_n+q'_n=q''_n}v_{p_1,\ldots,p_{n},q_1,\ldots,q_{n}}^{k,j}a_{p'_1,\ldots,p'_{n},q'_1,\ldots,q'_{n}}^{j,l}\prod_{s=1}^{n-1}\prod_{r=s+1}^n\l_{r,s}^{p_rp'_s+q_rq'_s+q_sp'_r-q_rp'_s}\\
		=&\left\{\begin{aligned}
			v_{p''_1,\ldots,p''_{n},q''_1,\ldots,q''_{n}}^{k,l},~1\leq k\leq r,\\
			0,~r+1\leq k\leq N.\\
		\end{aligned}\right.
	\end{aligned}
\end{equation}
	Then for all  $(p'''_1,\ldots,p'''_{n},q'''_1,\ldots,q'''_{n})\in\N^{2n}$ with $\sum_{j=1}^n\left(p'''_j+q'''_j\right)\geq1$, 
		\begin{align*}
			&c_{p'''_1,\ldots,p'''_{n},q'''_1,\ldots,q'''_{n}}^{k,l}\\
			=&\sum_{\a,\b=1}^N\sum_{p_1+p'_1+p''_1=p'''_1,\ldots,p_n+p'_n+p''_n=p'''_n\atop q_1+q'_1+q''_1=q'''_1,\ldots,q_n+q'_n+q''_n=q'''_n} v_{p_1,\ldots,p_{n},q_1,\ldots,q_{n}}^{k,\a}a_{p'_1,\ldots,p'_{n},q'_1,\ldots,q'_{n}}^{\a,\b}w_{p''_1,\ldots,p''_{n},q''_1,\ldots,q''_{n}}^{\b,l}\\
			&\times\prod_{s=1}^{n-1}\prod_{r=s+1}^n\l_{r,s}^{p_rp'_s+q_rq'_s+q_sp'_r-q_rp'_s+p'''_rp''_s+q'''_rq''_s+\left(q'''_s-q''_s\right)p''_r-\left(q'''_r-q''_r\right)p''_s}\\
			=&\sum_{\a,\b=1}^N\sum_{p_1+p'_1+p''_1=p'''_1,\ldots,p_n+p'_n+p''_n=p'''_n\atop q_1+q'_1+q''_1=q'''_1,\ldots,q_n+q'_n+q''_n=q'''_n} v_{p_1,\ldots,p_{n},q_1,\ldots,q_{n}}^{k,\a}a_{p'_1,\ldots,p'_{n},q'_1,\ldots,q'_{n}}^{\a,\b}\prod_{r=s+1}^n\l_{r,s}^{p_rp'_s+q_rq'_s+q_sp'_r-q_rp'_s}\\
			&\times w_{p''_1,\ldots,p''_{n},q''_1,\ldots,q''_{n}}^{\b,l}\prod_{s=1}^{n-1}\prod_{r=s+1}^n\l_{r,s}^{\left(p_r+p'_r+p''_r\right)p''_s+\left(q_r+q'_r+q''_r\right)q''_s+\left(q_s+q'_s\right)p''_r-\left(q_r+q'_r\right)p''_s}\\
			=&\left\{\begin{aligned}
				\sum_{\b=1}^N\sum_{p'_1+p''_1=p'''_1,\ldots,p'_n+p''_n=p'''_n\atop q'_1+q''_1=q'''_1,\ldots,q'_n+q''_n=q'''_n} v_{p'_1,\ldots,p'_{n},q'_1,\ldots,q'_{n}}^{k,\b}w_{p''_1,\ldots,p''_{n},q''_1,\ldots,q''_{n}}^{\b,l}\prod_{r=s+1}^n\l_{r,s}^{p''_rp''_s+q''_rq''_s+p'_rp''_s+q'_rq''_s+q'_sp''_r-q'_rp''_s}\\
				=0,~1\leq k\leq r,\\
				0,~r+1\leq k\leq N,\\
			\end{aligned}\right.
		\end{align*}
		which means
		\begin{equation}
			\mathcal{U}\mathcal{P}\mathcal{U}^*=\begin{pmatrix}
				I_r & 0\\
				0 & 0\\
			\end{pmatrix}.
		\end{equation}	
	Hence we have $K_0(\incren)=\Z$.
	\end{proof}
	For the case where $m$ is odd, $m=2n+1$, replace all $a_{p_1,\ldots,p_{n},q_1,\ldots,q_{n}}^{k,l}$ with series of $x$ in the proof of Theorem \ref{inp2n}, and $K_0(\C[x])=\Z$ guarantees that we can still make the assumption \eqref{asa}. Then by applying the same method as above, we can prove $K_0(\incron)=\Z$, and we finish the proof of Theorem \ref{m}.
	\begin{remark}
		We left the study of Morita equivalence of $m$-planes for future.
	\end{remark}

\end{document}